
\input amstex.tex
\documentstyle{amsppt}
\magnification1200
\hsize=12.5cm
\vsize=18cm
\hoffset=1cm
\voffset=2cm

\def\DJ{\leavevmode\setbox0=\hbox{D}\kern0pt\rlap
{\kern.04em\raise.188\ht0\hbox{-}}D}
\def\dj{\leavevmode
 \setbox0=\hbox{d}\kern0pt\rlap{\kern.215em\raise.46\ht0\hbox{-}}d}

\def\txt#1{{\textstyle{#1}}}
\baselineskip=13pt
\def\hf{{\textstyle{1\over2}}}

\def\d{{\,\roman d}}
\def\e{\varepsilon}

\def\f{\varphi}
\def\G{\Gamma}

\def\s{\sigma}
\def\t{\theta}
\def\={\;=\;}
\def\zx{\zeta(\hf+ix)}
\def\zt{\zeta(\hf+it)}

\def\L{\Cal L}

\def\R{\Re{\roman e}\,} 
\def\z{\zeta}

\def\t{\theta}
\def\hf{{\textstyle{1\over2}}}
\def\txt#1{{\textstyle{#1}}}
\def\f{\varphi}
\def\Z{{\Cal Z}}
\def\M{{\Cal M}}
\def\L{{\Cal L}}
\def\le{\leqslant}
\def\ge{\geqslant}
\font\tenmsb=msbm10
\font\sevenmsb=msbm7
\font\fivemsb=msbm5
\newfam\msbfam
\textfont\msbfam=\tenmsb
\scriptfont\msbfam=\sevenmsb
\scriptscriptfont\msbfam=\fivemsb
\def\Bbb#1{{\fam\msbfam #1}}

\def \NN {\Bbb N}
\def \CC {\Bbb C}
\def \RR {\Bbb R}

\font\ff=cmr8
\def\txt#1{{\textstyle{#1}}}
\baselineskip=13pt

\font\teneufm=eufm10
\font\seveneufm=eufm7
\font\fiveeufm=eufm5
\newfam\eufmfam
\textfont\eufmfam=\teneufm
\scriptfont\eufmfam=\seveneufm
\scriptscriptfont\eufmfam=\fiveeufm
\def\mathfrak#1{{\fam\eufmfam\relax#1}}

\font\tenmsb=msbm10
\font\sevenmsb=msbm7
\font\fivemsb=msbm5
\newfam\msbfam
     \textfont\msbfam=\tenmsb
      \scriptfont\msbfam=\sevenmsb
      \scriptscriptfont\msbfam=\fivemsb
\def\Bbb#1{{\fam\msbfam #1}}

\def \NN {\Bbb N}
\def \CC {\Bbb C}

\def \RR {\Bbb R}

  \def\rightheadline{{\hfil{\ff
Mellin transforms of Hardy's function}\hfil\tenrm\folio}}

  \def\leftheadline{{\tenrm\folio\hfil{\ff
 Aleksandar Ivi\'c }\hfil}}
  \def\emptyheadline{\hfil}
  \headline{\ifnum\pageno=1 \emptyheadline\else
  \ifodd\pageno \rightheadline \else \leftheadline\fi\fi}

\topmatter
\title
On some relations for Mellin transforms of Hardy's function
\endtitle
\author
 Aleksandar Ivi\'c
\endauthor
\address
Katedra Matematike RGF-a,
Universitet u Beogradu,  \DJ u\v sina 7,
11000 Beograd, Serbia.\bigskip
\endaddress
\subjclass 11 M 06
\endsubjclass
\email
{\tt
ivic\@rgf.bg.ac.rs,\enskip aivic\@matf.bg.ac.rs}
\endemail
\dedicatory
To Professor Akio Fujii on the occasion of his retirement
\enddedicatory
\abstract
Some relations involving the Mellin and Laplace
transforms of powers of the classical Hardy function
$$
Z(t) := \zt\bigl(\chi(\hf+it)\bigr)^{-1/2}, \quad \z(s) = \chi(s)\z(1-s)
$$
are obtained. In particular, we discuss some mean square identities and
their consequences.
\endabstract
\endtopmatter
\document
\head
1. Definition of Hardy's function
\endhead

\medskip
The familiar
Riemann zeta-function
$$
\z(s) = \sum\limits_{n=1}^\infty n^{-s}\qquad(s=\s +it, \s>1)
$$
admits analytic continuation to $\CC$, having only a simple pole at $s=1$.
A vast literature
exists on many aspects of zeta-function theory,
such as the distribution of its zeros and power
moments of $|\zt|$ (see e.g., the monographs [9], [10], [38], [39] and [42]).
It is within this framework that
the classical Hardy function (see e.g., [9]) $Z(t)\; (t\in\RR)$
plays an important r\^ole. It is defined as
$$
Z(t) := \zt\bigl(\chi(\hf+it)\bigr)^{-1/2},\leqno(1.1)
$$
where $\chi(s)$ comes from the well-known functional equation for $\z(s)$
(see e.g., [9, Chapter 1]),
namely $\z(s) = \chi(s)\z(1-s)$, so that
$$
\chi(s) = 2^s\pi^{s-1}\sin(\hf \pi s)\G(1-s),\quad
\chi(s)\chi(1-s)=1.
$$
It follows  that
$$
\overline{\chi(\hf + it)} = \chi(\hf-it)= \chi^{-1}(\hf+it),
$$
so that $Z(t)\in\RR$ when $t\in\RR$,
and $|Z(t)| =|\zt|$. Thus the zeros of $\z(s)$ on the ``critical line'' $\R s =1/2$
correspond to the real zeros of $Z(t)$, which makes $Z(t)$ an invaluable tool
in the study of the zeros of the zeta-function on the critical line.
Alternatively, if we use the symmetric form of the functional equation
for $\z(s)$, namely
$$
\pi^{-s/2}\z(s)\G(\hf s) = \pi^{-(1-s)/2}\z(1-s)\G(\hf(1-s)),
$$
then for $t\in \RR$ we obtain an equivalent form of (1.1). This is
$$
Z(t) = {\roman e}^{i\t(t)}\zt,\quad {\roman e}^{i\t(t)} := \pi^{-it/2}\frac
{\G({1\over4}+\hf it)}{|\G({1\over4}+\hf it)|},\quad \t(t) \in\RR.
$$
Hardy's original application of $Z(t)$ was to show that $\z(s)$ has infinitely
many zeros on the critical line $\R s =1/2$  (see e.g., E.C. Titchmarsh [42]).
Later A. Selberg (see [40] and [41]) obtained his famous result that a positive
proportion of zeros of $\z(s)$ lies on the critical line.
His method was later used and refined by many mathematicians.
The latest result is by S. Feng [5], who proved that at least 41.73\% of the zeros of $\z(s)$
are on the critical line and at least 40.75\% of  the zeros of $\z(s)$
are simple and on the critical line.

\smallskip
There is also an important connection of $Z(t)$ to the Riemann Hypothesis
(RH, that all complex zeros of $\z(s)$ have real parts 1/2).
The function $Z(t)$ has a negative local maximum
$-0.52625\ldots$ at $t = 2.47575\ldots\,$. This is the only known
occurrence of a negative local maximum, while no positive local
minimum is known. The so-called {\it Lehmer's phenomenon} (named after D.H. Lehmer,
who in his works [35], [36] made significant contributions to
the subject) is the fact that the graph of $Z(t)$ sometimes barely crosses the $t$--axis.
This means that the absolute value of the maximum or minimum
of $Z(t)$ between its two consecutive zeros is small.
The Lehmer phenomenon shows the delicacy
of the RH, and the possibility that a counterexample to the RH may be
found numerically. For should it happen that, for $t \ge t_0$,
$Z(t)$ attains a negative local maximum or a positive local minimum,
then the RH would be disproved. This assertion follows (see [18])
from the following proposition: If the RH is true, then the graph of
$Z'(t)/Z(t)$ is monotonically decreasing between the zeros  of $Z(t)$ for $t \ge t_0$.

\medskip
The main aim of this article is to discuss the
Mellin transforms of $Z^k(t)$ when
$k\in\NN$. The (modified) Mellin transform is introduced in Section 2, and
some of its properties are given. New mean square results for $\M_1(s)$ and
$\M_2(s)$ are presented in Section 3 and their proofs are given in Section 4.
Section 5 is devoted to the Laplace transforms of $Z^k(t)$, while some
other relations and open problems are stated in Section 6.
\medskip
\head
2. The modified  Mellin transform
\endhead
\medskip

First we recall that the Laplace and Mellin transforms of $f(x)$ are commonly defined as
($\s,\,t\in \Bbb R)$
$$
\L[f(x)] = \int_0^\infty f(x){\roman e}^{-sx}\d x \qquad(s = \s + it), \leqno(2.1)
$$
$$
\M[f(x)] =  F(s) := \int_0^\infty f(x)x^{s-1}\d x \qquad(s = \s + it),\leqno(2.2)
$$
provided that the integrals in question exist.
Mellin and Laplace transforms play an important r\^ole in Analytic Number Theory.
They can be viewed, by  a change of variable, as special cases of
Fourier transforms, and their properties can be deduced from the
general theory of Fourier transforms (see e.g., E.C. Titchmarsh [41]).

\smallskip
One of the basic properties of Mellin transforms is the  inversion formula
$$
\hf\{f(x+0) + f(x-0)\} = {1\over2\pi i}\int\limits_{(\s)}F(s)x^{-s}\d s =
{1\over2\pi i}\lim_{T\to\infty}\int\limits_{\s-iT}^{\s+iT}F(s)x^{-s}\d s,\leqno(2.3)
$$
where $\int_{(\s)}$ denotes integration over the line $\R s = \s$.
Formula (2.3) certainly holds  if $f(x)x^{\s-1} \in L^1(0,\infty)$, and $f(x)$
is of bounded variation on every finite $x$--interval.
Note that if $G(s)$ denotes the Mellin transform of $g(x)$  then,
assuming $f(x)$ and $g(x)$ to be real-valued, we formally have
$$
\eqalign{& {1\over2 \pi i}\int_{(\s)}F(s)\overline{G(s)}\d s
=\int_0^\infty g(x)\left({1\over2 \pi
i}\int_{(\s)}F(s)x^{\s-it-1}\d s\right)\d x\cr  &= \int_0^\infty
g(x)x^{2\s-1}\left({1\over2 \pi i}\int_{(\s)} F(s)x^{-s}\d s\right)\d x
=\int_0^\infty f(x)g(x)x^{2\s-1}\d x.\cr}
\leqno(2.4)
$$
The relation (2.4) is a form of Parseval's formula for Mellin
transforms, and it offers various possibilities  for mean square
bounds. A condition under which (2.4) holds is that $x^\s f(x)$ and $x^\s g(x)$
belong to $L^2((0,\infty),\, \d x/x)$, where as usual
$$
L^p(a,b) := \left\{\,f(x) \;\Biggl|\; \int_a^b|f(x)|^p\d x < \infty\;\right\}.
$$

\smallskip

Our main object of study will be the function
$$
{\Cal M}_k(s) :=
\int_1^\infty Z^k(x)x^{-s}\d x\qquad(k\in\NN),
\leqno(2.5)
$$
where  $\s = \R s = \s(k)$ is so large that the integral in (2.5) converges
absolutely, and henceforth $k\in\NN$ will be fixed. This function, which
appropriately can be called the {\it modified Mellin transform} of $Z^k(t)$,
was introduced and investigated in [19], [20].

The reasons why we  have defined  somewhat differently ${\Cal M}_k(s)$ from
the usual  Mellin transforms are the following: the lower limit of integration
$x=1$ dispenses with potential convergence problems at $x=0$,
while the appearance
of $x^{-s}$ instead of the familiar $x^{s-1}$ stresses the analogy with Dirichlet
series where one has a sum of $f(n)n^{-s}$ and not $f(n)n^{s-1}$. Note
that in the case
when $k = 2m$ is even, in view of $|Z(t)| = |\zt|$, we have
$$
\M_{2m}(s) = \int_1^\infty |\zx|^{2m}x^{-s}\d x.
$$
The special cases $m =1 $ and $m = 2$ were investigated in several works,
including [13], [15], [17], [21], [25], [26], [37] and  [38].
\medskip
The general modified Mellin transform $m[f(x)]$,
of which $\M_k(s)$ is a special case, is defined as
$$
F^*(s) \= m[f(x)] \= \int_1^\infty f(x)x^{-s}\d x\qquad(s = \s+it),
\leqno(2.6)
$$
provided that the integral in (2.6) converges.
 If ${\bar f}(x) = f(1/x)$ when $0 < x \le 1$ and
${\bar f}(x) = 0$ otherwise, then
$$
m[f(x)] \= \M\left[{1\over x}{\bar f}(x)\right],\leqno(2.7)
$$
so that the properties of $m[f(x)]$ can be deduced from the properties
of the ordinary Mellin transform $\M[f(x)]$ by the use of (2.7).
In this way the author [13] proved several lemmas involving the
 modified Mellin transform. We present three of them, which will
 be used in the sequel.

\medskip
LEMMA 1. {\it If $x^{-\s}f(x) \in L^1(1,\infty)$ and $f(x)$ is
continuous for $x > 1$, then}
$$
f(x) = {1\over2\pi i}\int_{(\s)}F^*(s)x^{s-1}\d s,\quad F^*(s) \=
m[f(x)]. \leqno(2.8)
$$

\smallskip
The inversion formula (2.8) is the analogue of (2.3) for modified
Mellin transforms.
The next lemma is the analogue of (2.4) for  modified Mellin transforms.


\smallskip
LEMMA 2. {\it If $F^*(s) = m[f(x)],\,G^*(s) = m[g(x)]$, and $f(x), g(x)$ are
real-valued, continuous functions for $x > 1$, such that}
$$
x^{{1\over2}-\s}f(x) \in L^2(1,\infty),\quad x^{{1\over2}-\s}g(x)
\in L^2(1,\infty),
$$
{\it then}
$$
\int_1^\infty f(x)g(x)x^{1-2\s}\d x =  {1\over2\pi i}\int_{(\s)}
F^*(s)\overline{G^*(s)}\d s. \leqno(2.9)
$$

\medskip
LEMMA 3. {\it Suppose that $g(x)$ is a real-valued,
integrable function on $[a,b]$, a subinterval
of $[2,\,\infty)$, which is not necessarily finite. Then}
$$
\int\limits_0^{T}\left|\int\limits_a^b g(x)x^{-s}\d x\right|^2\d t
\le 2\pi\int\limits_a^b g^2(x)x^{1-2\s}\d x \quad(s = \s+it\,,T > 0,\,a<b).
\leqno(2.10)
$$
\medskip

\medskip
\head
3. Mean square identities with $\M_k(s)$
\endhead
To formulate our results on
some mean square identities with $\M_k(s)$, we need
the definition of the function $E_k(T)$, the error
term in the asymptotic formula for the $2k$-th moment
of $|\zt|$. Namely
for any fixed $k \in\NN$, we expect (since the lower bound of integration
in (2.5) is unity, it is convenient to have it also in the integral
in (3.1))
$$
\int_1^T|\zt|^{2k}\d t = TP_{k^2}(\log T) + E_k(T)\leqno(3.1)
$$
to hold, where it is generally assumed that
$$
P_{k^2}(y) = \sum_{j=0}^{k^2}a_{j,k}y^j\leqno(3.2)
$$
is a polynomial in $y$ of degree $k^2$ (the integral in (3.1) is
unconditionally $\gg_k T\log^{k^2}T$; see e.g., [9, Chapter 9]
and [39]). The function $E_k(T)$ is to be considered as the error
term in (3.1), namely one supposes that
$$
E_k(T) \= o(T)\qquad(T \to \infty).\leqno(3.3)
$$
So far the formulas (3.1)--(3.3) are known to hold only for $k = 1$ and $k = 2$
(see [9] and [10] for a detailed account). For higher moments
one has the bound
$$
\int_1^T|\zt|^{2k}\d t \;\ll\; T^{{1\over4}(k+2)}\log^{C_k}T\qquad(2\le k\le 6),
\leqno(3.4)
$$
and more complicated bounds for $k>6$. A comprehensive account is
to be found in Chapter 8 of [9].
As usual, $f \ll g$ (same as $f = O(g)$), means that $|f(x)| \le Cg(x)$
for some $C>0$ and $x\ge x_0$, and $f \ll_\e g$ means that the
$\ll$-constant depends on $\e$. Therefore in view of the
existing knowledge on the higher moments of $|\zt|$,
embodied essentially in (3.4), at present
the really important cases of (3.1) are $k = 1$ and $k = 2$. Plausible
heuristic arguments for the values of the coefficients $a_{j,k}$
in the general case were
given by Conrey et al. [4], by using methods from Random Matrix Theory
(see also Keating--Snaith [32]).

\smallskip
For $k=1$ the relation (3.1) becomes (see Chapter 15 of [9])
the well-known mean square formula
$$
\int_1^T|\zt|^2\d t = T\Bigl(\log{T\over2\pi} + 2\gamma -1\Bigr) + E(T),\leqno(3.5)
$$
where one usually writes $E(T) \equiv E_1(T)$, and
$\gamma = -\G'(1) = 0.57721 56649 \ldots$ is
Euler's constant. It is known that
$E_1(T) \ll_\e T^{131/416+\e}$ (see Huxley-Ivi\'c [8] and N. Watt [43])
 and on the other hand we have $E_1(T) = \Omega_\pm(T^{1/4})$.
 Also we have (see Hafner-Ivi\'c [6], [7])
$$
\int_1^T E(t)\d t = \pi T + G(T),\quad G(T) = O(T^{3/4}),\quad G(T) = \Omega_\pm(T^{3/4}).
\leqno(3.6)
$$

\smallskip
For $k=2$ we have
$$
E_2(T) \;=\; \int_1^T|\zt|^4\d t - TP_4(\log T).\leqno(3.7)
$$
Note that $P_4(x)$ is a polynomial
of degree four in $x$ with leading coefficient $1/(2\pi^2)$
(see e.g., [9, Chapter 4]). Its coefficients
were evaluated independently by J.B. Conrey [3] and the author [11].
We have $E_2(T) \ll T^{2/3}\log^9T$ and $E_2(T) = \Omega_\pm(T^{1/2})$
(see [23] and [38]). It is conjectured that $E_1(T) \ll_\e T^{1/4+\e}$
and that $E_2(T) \ll_\e T^{1/2+\e}$.

\smallskip
Here and later $\e\;(>0)$ will denote arbitrarily small constants,
not necessarily the same ones at each occurrence. As usual,
$f = \Omega_\pm(g)$ means that $\limsup f/g > 0$ and $\liminf f/g < 0$.
We also have (see [10])
$$
\int_1^T E^2(t) \d t = DT^{3/2} + O(T\log^4T)\leqno(3.8)
$$
with $D = 2(2\pi)^{-1/2}\z^4(3/2)/(3\z(3))$, and (see [23] and [38])
$$
\int_1^{T} E^2_2(t) \d t   \ll T^2\log^{22}T.\leqno(3.9)
$$

\smallskip
With this notation we can formulate our mean square results for $\M_k(s)$.

\smallskip
THEOREM 1. {\it For $\s>1$ we have}
$$
\eqalign{
&{1\over\pi}\int_0^\infty|\M_1(\s+it)|^2\d t = {2\gamma-\log(2\pi)\over2\s-2}
+ {1\over(2\s-2)^2} + \Bigl(2\gamma-1-\log(2\pi)\Bigr) \cr&
+ 2\pi\s+ 2\s(2\s-1)\int_1^\infty G(x)x^{-1-2\s}\d x,\cr}\leqno(3.10)
$$
{\it where the integral on the right-hand side of} \,(3.10) {\it converges
absolutely for $\s > 3/8$.}

\medskip
THEOREM 2. {\it For $\s>1$ we have}
$$
{1\over\pi}\int_0^\infty|\M_2(\s+it)|^2\d t = \sum_{j=0}^5\,{c_j\over(\s-1)^j} +
(2\s-1)\int_1^\infty E_2(x)x^{-2\s}\d x,\leqno(3.11)
$$
{\it where the constants $c_j$ can be evaluated explicitly, and
the integral on the right-hand side of} (3.11) {\it converges
absolutely for $\s > 3/4$.}

 \medskip
 {\bf Corollary}. We have
$$
\eqalign{
&\lim_{\s\to1+0}\Bigl\{{1\over\pi}\int_0^\infty\M_1(\s+it)|^2\d t
 - {2\gamma-\log(2\pi)\over2\s-2}
- {1\over(2\s-2)^2}\Bigr\}\cr& = 2\pi + 2\gamma-1-\log(2\pi)+
2\int_1^\infty G(x)x^{-3}\d x\cr}
$$
and
$$
\lim_{\s\to1+0}\Bigl\{{1\over\pi}\int_0^\infty\M_2(\s+it)|^2\d t -
\sum_{j=0}^5\,{c_j\over(\s-1)^j}\Bigr\}
= \int_1^\infty E_2(x)x^{-2}\d x.
$$

\smallskip
{\bf Remark}. Theorem 1 could have been formulated, analogously to (3.11), as
$$
{1\over\pi}\int_0^\infty|\M_1(\s+it)|^2\d t =
C_0 + {C_1\over\s-1} + {C_2\over(\s-1)^2} + (2\s-1)\int_1^\infty E(x)x^{-2\s}\d x,
$$
where the integral on the right-hand side converges absolutely for $\s>5/8$.
This follows by the Cauchy-Schwarz inequality for integrals from (3.8),
but (3.10) is more precise.

\medskip
For $k\ge3$ let us define
$$
\t_k \;:=\; \inf\Bigl\{\;a_k \;:\; \int_{-\infty}^\infty
|\M_k(\s+it)|^2\d t < \infty\;\;\roman{for}\;
\;\s>a_k\Bigr\}.
$$
It is clear that $\t_k$ always exists (e.g., since $\zt \ll |t|^{1/6}$)
and that it has an intrinsic connection with power moments of $|\zt|$.
Of course, the above definition makes sense for $k=1,2$ as well,
 but in these cases we have much more precise
information in view of Theorem 1 and Theorem 2.

\medskip
THEOREM 3. {\it We have}
$$
\t_k\ge 1\quad(\forall k),\quad \t_k \le {3\over4} + {k\over8}\quad(3\le k\le 6),\leqno(3.12)
$$
{\it and}
$$
\int_1^T|\zt|^{2k}\d t \;\ll_\e\; T^{2\t_k-1+\e}\qquad(k\ge3).\leqno(3.13)
$$
\medskip
\head
4. Proofs of the theorems
\endhead
\medskip
We begin with (2.9) of Lemma 2, which in case when
$$
f(x) = g(x) = Z^k(x)
$$
reduces to
$$
\eqalign{
\int_1^\infty |\zx|^{2k}x^{1-2\s}\d x &= {1\over2\pi}\int_{-\infty}^\infty
|\M_k(\s+it)|^2\d t\cr&
= {1\over\pi}\int_{0}^\infty
|\M_k(\s+it)|^2\d t,\cr}\leqno(4.1)
$$
since $\overline{\M_k(s)} = \M_k({\bar s})$.
To evaluate the left-hand side of (4.1) note that differentiation of
(3.1) yields
$$
|\zt|^{2k} = P_{k^2}(\log t) + P'_{k^2}(\log t) + E'_k(t),
$$
with $P_{k^2}$ given by (3.2). Hence, initially for $\R\s \ge \s_1(k)$, we have

$$\eqalign{&
\int_1^\infty |\zx|^{2k}x^{1-2\s}\d x =
\int_1^\infty x^{1-2\s}\d\Bigl\{xP_{k^2}(\log x) + E_k(x)\Bigr\}\cr& =
\int\limits_1^\infty (P_{k^2}(\log x) + P'_{k^2}(\log x))x^{1-2\s}\d x -
E_k(1) + (2\s-1)\int\limits_1^\infty E_k(x)x^{-2\s}\d x.\cr}\leqno(4.2)
$$
But for $\R\s > 1$ change of variable $\log x = t$ gives
$$
\eqalign{& \int_1^\infty (P_{k^2}(\log x) + P'_{k^2}(\log x))x^{1-2\s}\d x\cr&
 = \int_1^\infty \left\{\sum_{j=0}^{k^2}a_{j,k}\log^jx +
\sum_{j=0}^{k^2-1}(j+1)a_{j+1,k}\log^jx\right\}x^{1-2\s}\d x\cr& =
\int_0^\infty \left\{\sum_{j=0}^{k^2}a_{j,k}t^j +
\sum_{j=0}^{k^2-1}(j+1)a_{j+1,k}t^j\right\}{\roman e}^{-(2\s-2)t}\d t
\cr& = {a_{k^2,k}(k^2)!\over(2\s-2)^{k^2+1}} +
\sum_{j=0}^{k^2-1}(a_{j,k}j! + a_{j+1,k}(j+1)!)(2\s-2)^{-j-1}.
\cr}\leqno(4.3)
$$

When $k=1$ we have, by (3.5),
$$
P_1(y) = y + 2\gamma-1-\log(2\pi),
$$
hence for $\R\s>1$
$$
\eqalign{&
\int_1^\infty\left(P_1(\log x)+P'_1(\log x)\right)x^{1-2\s}\d x
\cr&
= \int_1^\infty(\log x + 2\gamma-\log(2\pi))x^{1-2\s}\d x
\cr&=
{2\gamma-\log(2\pi)\over2\s-2} + {1\over(2\s-2)^2}.\cr}\leqno(4.4)
$$
Next note that
$$
E_1(1) \equiv E(1) = -P_1(0) = \log(2\pi) +1 -2\gamma.\leqno(4.5)
$$
Finally integration by parts yields, on using (3.6),
$$
\eqalign{&
\int_1^\infty E_1(x)x^{-2\s}\d x\cr&
= \int_1^x E(y)\d y\cdot x^{-2\s}\Bigl|_1^\infty + 2\s\int_1^\infty
\int_1^x E(y)\d y\cdot x^{-1-2\s}\d x\cr&
= 2\s\int_1^\infty (\pi x + G(x))x^{-1-2\s}\d x\cr&
= {2\pi\s\over 2\s-1}+2\s\int_1^\infty G(x)x^{-1-2\s}\d x,\cr}\leqno(4.6)
$$
and the last integral converges absolutely for $\s>3/8$ in view of
the $O$-bound in (3.6).
The assertion of Theorem 1 follows then from (4.1)--(4.6).

\medskip
For $k=2$ write
$$
P_4(y) = \sum_{j=0}^4a_{j,4}y^j = \sum_{j=0}^4A_jy^j, A_4 = 1/(2\pi^2), \leqno(4.7)
$$
so that $E_2(1) = - P_4(0) = - A_0$ by (3.7).
From (4.2) and (4.7) we infer that, for $\s>1$,
$$
\int_1^\infty (P_{4}(\log x) + P'_{4}(\log x))x^{1-2\s}\d x
= \sum_{j=1}^5 \,{B_j\over(\s-1)^j}
$$
with
$$
B_5 = {3\over8\pi^2}, \;B_j = A_{j-1}(j-1)!2^{-j} + A_jj!2^{-j}\quad(j=1,2,3,4).
$$
This clearly gives (3.11) of Theorem 2 with
$$
c_0 = A_0, \;c_j = A_{j-1}(j-1)!2^{-j} + A_jj!2^{-j}\;\;(j=1,2,3,4),
\; c_5 = {3\over8\pi^2}.
$$
The integral on the right-hand side of (3.11) converges absolutely for $\s>3/4$,
since by the Cauchy-Schwarz inequality for integrals and (3.2) we have
$$
\int_X^{2X}E_2(x)x^{-2\s}\d x \le {\left\{\int_X^{2X}E_2^2(x)\d x\int_X^{2X}x^{-4\s}\d x
\right\}}^{1/2} \ll_\e X^{3/2-2\s+\e} \le X^{-\e}
$$
for $\s > 3/4 + \e$.

\medskip
To prove Theorem 3 recall that it was stated after (3.1) that unconditionally
$$
\int_1^T|\zt|^{2k}\d t \;\gg_k\; T(\log T)^{k^2}.
$$
This implies that $\M_k(s)$ diverges for $s=1$, hence $\t_k\ge1$ must hold. On the other
hand we use (3.4) to deduce that, for $3\le k\le6$,
$$
\int_X^{2X}|\zx|^{2k}x^{1-2\s}\d x \;\ll\; X^{1-2\s}X^{(k+2)/4}(\log X)^{C_k}
$$
and $1 - 2\s + (k+2)/4 < 0$ for $\s > 3/4 + k/8$. This means that,
 in this range for $\s$, the integral
$$
\int_1^\infty |\zx|^{2k}x^{1-2\s}\d x
$$
converges. But then (4.1) holds and hence $\t_k\le 3/4 + k/8$ for $3\le k\le 6$,
proving Theorem 3. The conjecture $\t_k = 1\,(\forall k)$ is clearly equivalent to the
Lindel\"of hypothesis that $\zt \ll_\e |t|^\e$ (see (7.2) of [9] and Theorem 13.2 of
[42]).

\medskip
\head
5. The Laplace transform of $Z^k(t)$
\endhead
Let
$$
\L_k(s) := \int_1^\infty Z^k(x)\,{\roman e}^{-sx}\d x
\qquad(\s = \R s >0,\;k\in\NN)\leqno(5.1)
$$
denote the (modified) Laplace transform of $Z^k(x)$. Analogously to the modified
Mellin transform (2.5), this differs from the standard
definition of the Laplace transform in the lower bound of integration which is in
(5.1) unity, and not zero like in (2.1). This is convenient because, for $c>0$ and
$\s \ge \s_0(k)$, we obtain
$$
\eqalign{
\L_k(s) &= {1\over2\pi i}\int_{(c)}\G(w)\left(\int_1^\infty (sx)^{-w}Z^k(x)\d x
\right)\d w\cr&
= {1\over2\pi i}\int_{(c)}\G(w)s^{-w}\M_k(w)\d w.
\cr}\leqno(5.2)
$$
Here we used the well-known Mellin inversion formula (see e.g., the Appendix of [8])
$$
{\roman e}^{-x}= {1\over2\pi i}\int_{(c)}\G(w)x^{-w}\d w\qquad(\R x >0,\,c>0).
$$
Therefore by the inversion formula for modified Mellin transforms (see (2.8)
of Lemma 1) one has
$$
\G(s)\M_k(s) \;=\;\int_1^\infty \L_k\Bigl({1\over x}\Bigr) x^{-1-s}\d x
\qquad(\s \ge \s_0(k)).\leqno(5.3)
$$
This relation was used by the author in [17]  to show that $\M_2(s)$ has
meromorphic continuation to $\CC$. Its poles are $s=1$ of order two, and simple
poles at $s = -1, -3, -5, \ldots\,$. This result was independently proved also
by M. Lukkarinen [37]. In [19] the author proved that
$\M_1(s)$ has analytic continuation which is regular for $\R s >0$,
and M. Jutila [30] showed that $\M_1(s)$ is even entire.
The functions $\M_3(s)$ and $\M_4(s)$ were investigated in [19],[20],[21] and [27].

\smallskip
From (2.9) of Lemma 2 and  (5.3) we obtain then
$$
\int_1^\infty \L_k^2\Bigl({1\over x}\Bigr)x^{-1-2\s}\d x =
{1\over\pi}\int_0^\infty |\G(\s+it)|^2|\M_k(\s+it)|^2\d t
\quad(\s \ge \s_2(k)\;(>0)\,).
\leqno(5.4)
$$
Since, for fixed $\s>0$, we have
$$
|\G(\s+it)|^2\;\sim\;2\pi t^{2\s-1}{\roman e}^{-\pi t}\qquad(t\to+\infty)
$$
by Stirling's formula for the gamma-function, this means that the integrals
$$
\int_1^\infty \L_k^2\Bigl({1\over x}\Bigr)x^{-1-2\s}\d x,
\quad
\int_1^\infty t^{2\s-1}{\roman e}^{-\pi t}|\M_k(\s+it)|^2\d t
$$
both converge or both diverge for a given $\s\,(>0)$.

\medskip
One can, of course, consider  also the classical Laplace transform
$$
L_k(s) \;:=\; \int_0^\infty Z^k(x){\roman e}^{-sx}\d x\qquad(k \in
\NN,\,\R s > 0).
$$
E.C. Titchmarsh's well-known monograph [42, Chapter 7] provides a
discussion of $L_{2m}(s)$ when $s = \s$ is real and $\s \to 0+$,
especially detailed in the cases $m=1$ and $m=2$. Indeed, a
classical result of H. Kober [33] says that, as $\s \to 0+$,
$$
L_2(2\s) = {\gamma-\log(4\pi\s)\over2\sin\s} +
 \sum_{n=0}^Nc_n\s^n + O_N(\s^{N+1})
\leqno(5.5)
$$
for any given integer $N \ge 1$, where the $c_n$'s are effectively
computable constants and $\gamma$ is
Euler's constant. For complex values of $s$ the function $L_4(s)$
was studied by F.V. Atkinson [1], and more recently by M. Jutila
[27], who noted that Atkinson's argument gives
$$
L_2(s) = -i{\roman e}^{{1\over2}is}\bigl(\log(2\pi)-\gamma + ({\pi\over2}-s)i\bigr) +
2\pi {\roman e}^{-{1\over2}is}\sum_{n=1}^\infty d(n)\exp(-2\pi in{\roman e}^{-is}) +
\lambda_1(s)
$$
in the strip $0 < \R s < \pi$, where the function $\lambda_1(s)$
is holomorphic in the strip $|\R s| < \pi$. Moreover, in any strip
$|\R s| \le \theta$ with $0 < \theta < \pi$, we have
$$
\lambda_1(s) \;\ll_\theta \;(|s|+1)^{-1}.
$$

For $L_4(\s)$ F.V. Atkinson [2] obtained the asymptotic formula
$$
L_4(\s) = {1\over\s}\left(A\log^4{1\over\s} + B\log^3{1\over\s} +
C\log^2{1\over\s} + D\log {1\over\s} + E\right) +
\lambda_4(\s),\leqno(5.6)
$$
where $\s \to 0+$,
$$
A = {1\over2\pi^2},\,B =\pi^{-2}(2\log(2\pi) - 6\gamma +
24\z'(2)\pi^{-2})
$$
and
$$
 \lambda_4(\s) \;\ll_\e\;\left({1\over\s}\right)^{{13\over14}+\e}.
 \leqno(5.7)
$$
He also indicated how, by the use of estimates for Kloosterman
sums, one can improve the exponent ${13\over14}$ in (5.7) to
${8\over9}$. This is of historical interest, since it is one of
the first instances of an application of Kloosterman sums to analytic
number theory. Atkinson in fact showed that ($\s = \R s > 0$)
$$
L_4(s) \;=\;4\pi {\roman e}^{-{1\over2}s}\sum_{n=1}^\infty d_4(n) K_0(4\pi
i\sqrt{n}{\roman e}^{-{1\over2}s}) + \phi(s),\leqno(5.8)
$$
where $d_4(n)$ is the divisor function generated by $\z^4(s)$,
$K_0$ is the familiar Bessel function, and the series in (5.8) as well as
the function $\phi(s)$ are both analytic in the region $|s| <
\pi$.

\medskip
Note that
the author [11] applied a result on the fourth moment of $|\zt|$,
obtained jointly with Y. Motohashi [22]--[24] (see also [38]), to
establish that
$$
\lambda_4(\s) \;\ll\; \s^{-1/2}\qquad(\s\to 0+).
$$
This is essentially best possible, as shown by the author in [14],
who obtained a refinement of (5.6) by means of the spectral theory
of the non-Euclidean Laplacian (see Y.
Motohashi's monograph [38] for a comprehensive account).

\medskip
For $k \ge 5$ not much is known about $L_k(s)$, even when $s = \s
\to 0+$. This is not surprising, and is analogous to the
situation with $\M_k(s)$, since not much is known (cf. (3.4)) about
upper bounds for the $k$-th moment of $|\zt|$ when $k\ge5$.

\medskip
\head
6. Further discussion and some open problems
\endhead
There is a natural connection between $\M_k(s)$ and the moments
of $|\zt|$. For example, the author [19] proved that
$$
\int_T^{2T}|\zt|^6\d t \ll_\e T^{2\s-1}
\int_1^{T^{1+\e}}|{\Cal M}_3(\s+it)|^2\d t
+ T^{1+\e}\quad(\hf < \s \le 1), \leqno(6.1)
$$
provided that  ${\Cal M}_3(s)$ can be continued analytically to $\R s \ge\s$
(and that is the catch!).
Heuristically, we should be able to have $\s = 3/4+\e$, and then the integral on the
right-hand side of (6.1) should be $\ll_\e T^{1/2+\e}$,
giving the bound $O_\e(T^{1+\e})$, which is a weak form of the sixth moment.
Note that  (see [13, eq. (4.7)]) for the eighth moment we have
(since $\Z_2(s) \equiv \M_4(s)$)
$$
\int_T^{2T}|\zt|^8\d t \ll_\e
T^{2\s-1}\int_1^{T^{1+\e}}|\M_4(\s+it)|^2\d t + T^{1+\e}\quad(\hf < \s \le1),
\leqno(6.2)
$$
and an analogue of (6.1) and (6.2) holds also for the mean square and
fourth power of $|\zt|$. In these cases, however, the results are not
of particular interest, since we have precise information which has been
obtained by other methods. Note that $\M_4(s)$, unlike $\M_3(s)$,
is known to possess analytic continuation
to the region $\s > \hf$, where it is regular except for a pole of order
five at $s=1$ (see Y. Motohashi [37]).

\smallskip
The (nontrivial) bounds for the sixth moment of $|\zt|$ are
intricately connected to the problem of the analytic continuation of
$\M_3(s)$ to the region $\s \le 1$. This, in turn, depends on the asymptotic
evaluation of the integral
$$
F_k(T) \;:=\;\int_1^T Z^k(t)\d t\leqno(6.3)
$$
when $k=3$. The author in [16] proved that
$$
F_1(T) = \int_1^T Z(t)\d t = O_\e(T^{1/4+\e}),
$$
which was improved to $F_1(T) = O(T^{1/4})$ by M. Korolev [34],
who also proved that $F(T) = \Omega_\pm(T^{1/4})$.
M. Jutila [28], [31] gave a different proof of the same
results by establishing precise formulas for $F_1(T)$.
In [19] it was proved that, for $k = 1,2,3,4$, we have
$$
\eqalign{
&F_k(2T) - F_k(T) = \int\limits_T^{2T}Z^k(t)\d t\cr&= 2\pi\sqrt{2\over k}
\sum_{({T\over 2\pi})^{k/2}\le n\le ({T\over\pi})^{k/2}}
d_k(n)n^{-{1\over2}+{1\over k}}\cos\bigl(k\pi n^{2\over k}+{\txt{1\over8}}(k-2)\pi\bigr)
+ \cr& +\ldots +O_\e(T^{k/4+\e}),\cr}\leqno(6.4)
$$
where $+\ldots$ denotes terms similar to the one on the right-hand side of (6.4),
 with the similar cosine term, but of a lower order of magnitude. It was also indicated
 that actually the terms standing for $+\ldots+$ may be omitted. Here $d_k(n)$ is the
 divisor function generated by $\z^k(s)$ (so that $d_1(n)\equiv 1, d_2(n) \equiv
 d(n) = \sum_{\delta|n}1$). The  interesting case of (6.4)
 is  $k=3$ (since $k=1$ is solved,
 and for $k=2,4$ we have the well-known even moments), when the
 exponential sum in (6.4) can be estimated. This in turn furnishes the following result
 on $\M_3(s)$ (see [19, Theorem 5]): we have
$$
{\Cal M}_3(s) = \int_1^\infty Z^3(x)x^{-s}\d x = V_1(s) + V_2(s),
$$
say, where $V_2(s)$ is  regular for $\s > 3/4$,
and for $\s>1$ the function
$$
V_1(s) = (2\pi)^{1-s}\sqrt{2\over 3}\sum_{n=1}^\infty
d_3(n)n^{-{1\over6}-{2s\over3}}\cos\bigl(3\pi n^{2\over 3}+{\txt{1\over8}}\pi\bigr)
$$
is regular. In connection with this, one may naturally pose the following
problems (see [19], [20] for the first one):

\smallskip
1. Does there exist a constant $0< c_3<1$ such that
$$
F_3(T) = O(T^{c_3})?\leqno(6.5)
$$
Note that $c_3 = 1+\e$ is trivial, since by the Cauchy-Schwarz inequality for
integrals we easily obtain a better result, namely
$$
\left|\,\int_0^{T}Z^3(t)\d t\,\right| \le \left(\int_0^{T}|\zt|^2\d t
\int_0^{T}|\zt|^4\d t\right)^{1/2} \ll T(\log T)^{5/2}
$$
on using the well-known elementary bounds (see e.g., [9])
$$
\int_0^{T}|\zt|^2\d t \;\ll\; T\log T,\quad \int_0^{T}|\zt|^4\d t \;\ll\; T\log^4T.
$$

\smallskip
2. What is the least lower bound for $c_3$? It seems reasonable to conjecture that
(6.5) holds with $c_3 = 3/4+\e$ but does not hold for $c_3 < 3/4$, but so far
no positive lower bounds for $c_3$ are known.

\smallskip
3. Does there exist a constant $0< c_5\le9/8$ such that
$$
F_5(T) = O(T^{c_5})?\leqno(6.6)
$$
Namely (6.6) holds certainly with $c_5 = 9/8+\e$ in view of [9, eq. (8.57)], and any value
$c_5 \le 9/8$ would be non-trivial and very interesting.

\smallskip
4.  What is the least lower bound for $c_5$? Is it perhaps $c_5 = 1$, or is there
enough cancelation in the terms of $Z^5(t)$ to produce a smaller exponent than
unity? Naturally, similar questions could be asked for any   odd $k>3$, but
they are quite difficult in the general case.

\medskip
Note that by (2.10) of Lemma 2 we have, taking $f(x) = Z(x), g(x) = Z^4(x)$,
$$
\int_1^\infty Z^5(x)x^{1-2\s}\d x = {1\over2\pi i}\int_{(\s)}\M_1(s)\overline{\M_4(s)}
\d s\leqno(6.7)
$$
for $\s\,(>1)$ sufficiently large. Since integration by parts shows that (cf. (6.3))
$$
\int_1^\infty Z^5(x)x^{1-2\s}\d x = (2\s-1)\int_1^\infty F_5(x)x^{-2\s}\d x,
\leqno(6.8)
$$
then if $F_5(x) \ll_\e x^{c_5+\e}$, this implies that the left-hand side
of (6.8) converges for $\s > \hf(1+c_5)$, and in this range the integral on the
right-hand side of (6.7) converges as well. The integral on the right-hand side
of (6.7) can be dealt with by several techniques. One way is to use the Cauchy-Schwarz
inequality, the defining relations for $\M_k$ and Lemma 3. However, so far I have not
been able to improve on $c_5 \le 9/8+\e.$

\medskip

It may be also mentioned that in [19] it was shown that,
if $k = 1,2,3,4$ and $c>1$ is fixed, then for $U\gg x$
 and $\e>0$ sufficiently small,
 $$
 Z^k(x) = {1\over2\pi i}\int_{c-iU}^{c+iU}x^{s-1}\M_k(s)\d s + O_{\e,k}(x^{c-1}U^{-\e/2}).
 \leqno(6.9)
 $$
One may use (6.9) for various estimates involving moments of $Z(t)$. For example,
take $U = C_1X, X/2\le x\le 5X/2$, let $\f(x)\;(\ge0)$ be a smooth function
supported in $[X/2,\,5X/2]$ such that $\f(x) = 1$ when $X\le x\le2X$ and
$\f^{(r)}(x) \ll_r X^{-r}$ for $r\in\NN$. Then (6.9) yields
$$
\int\limits_{X/2}^{5X/2}\f(x)Z^k(x)\d x =
{1\over2\pi i}\int\limits_{c-iC_1X}^{c+iC_1X}\M_k(s)
\left(\int\limits_{X/2}^{5X/2}\f(x)x^{s-1}\d x\right)\d s + O_\e(X^{c-\e/2}).
$$
However, for any $r \,(\in\NN)$ repeated integration by parts yields
$$
\int\limits_{X/2}^{5X/2}\f(x)x^{s-1}\d x = (-1)^r\int\limits_{X/2}^{5X/2}\f^{(r)}(x)
{x^{s+r-1}\over s(s+1)\cdots (s+r-1)}\d x \ll_r {x^{\s-1}\over|t|^r} \ll x^{-A}
$$
for any given large $A>0$  provided that $|t|\ge X^\e$ if $r = [(A+\s)/\e].$
This gives, on writing $\s\,(>1)$ in place of $c$, for $k = 1,2,3,4$,
$$
\int\limits_{X/2}^{5X/2}\f(x)Z^k(x)\d x \ll_\e X^\s\max_{|t|\le X^\e}|\M_k(\s+it)|
 + X^{\s-\e/2}.
\leqno(6.10)
$$
The bound (6.10) shows essentially that the integral on the left-hand side is
bounded by $ X^{\s+\e}$,
if $\M_k(s)$ can be continued analytically to $\R s \ge\s$. This is in fact
another way of seeing how the power moments of $Z^k(t)$ and $\M_k(s)$ are connected.
Probably (6.10) holds for $k>4$ as well.
\vfill\eject
\topglue1cm
\Refs
\vskip1cm
\smallskip

\item{[1]} F.V. Atkinson, The mean value of the zeta-function on
the critical line, Quart. J. Math. Oxford {\bf 10}(1939), 122-128.

\item {[2]} F.V. Atkinson, The mean value of the zeta-function on
the critical line, Proc. London Math. Soc. {\bf 47}(1941), 174-200.

\item{[3]} J.B. Conrey, A note on the fourth power moment of the Riemann zeta-function,
B. C. Berndt (ed.) et al., in ``Analytic number theory. Vol. 1.
Proceedings of a conference in honor of Heini Halberstam, Urbana, 1995'',
 Birkhäuser, Prog. Math. {\bf138}(1996), 225-230.

\item{[4]} J.B. Conrey, D.W. Farmer, J.P. Keating, M.O. Rubinstein
and N.C. Snaith, Integral moments of $L$-functions,
Proc. London Math. Soc. (3) {\bf91}(2005), 33-104.

\item{[5]}  S. Feng, Zeros of the Riemann zeta-function
on the critical line, preprint available at \break arXiv:1003.0059.

\item{[6]}  J.L. Hafner  and A. Ivi\'c, On some mean value results for the
  Riemann zeta-function,  Proceedings International Number Theory Conference
  Qu\'ebec 1987, Walter de Gruyter and Co., 1989, Berlin - New York, 348-358.

\item{[7]} J.L. Hafner  and A. Ivi\'c, On the mean square of the Riemann zeta-function
on the critical line,  J. Number Theory   {\bf 32}(1989), 151-191.

\item{[8]} M.N. Huxley and A. Ivi\'c,
Subconvexity for the Riemann zeta-function and the divisor problem,
Bulletin CXXXIV de l'Acad\'emie Serbe des Sciences et des
Arts - 2007, Classe des Sciences math\'ematiques et naturelles,
Sciences math\'ematiques No. {\bf32}, pp. 13-32.

\item{[9]} A. Ivi\'c, The Riemann zeta-function, John Wiley \&
Sons, New York 1985 (2nd edition.  Dover, Mineola, New York, 2003).

\item {[10]} A. Ivi\'c,  Mean values of the Riemann zeta-function,
LN's {\bf 82},  Tata Inst. of Fundamental Research,
Bombay,  1991 (distr. by Springer Verlag, Berlin etc.).

\item{[11]}  A. Ivi\'c,  On the fourth moment of the Riemann
zeta-function,  Publs. Inst. Math. (Belgrade) {\bf 57(71)}
(1995), 101-110.

\item{[12]} A. Ivi\'c, The Laplace transform of the fourth
moment of the zeta-function,
Univ. Beograd. Publ. Elektrotehn. Fak. Ser. Mat. {\bf11}(2000), 41-48.

\item{[13]} A. Ivi\'c, On some conjectures and results for the
Riemann zeta-function, Acta. Arith. {\bf99}(2001), 115-145.

\item{[14]} A. Ivi\'c, The Laplace transform of the fourth moment of
of the zeta-function,  Univ. Beograd. Publ. Elektrotehn. Fak. Ser. Mat.
{\bf11}(2000), 41-48.

\item{ [15]}  A. Ivi\'c, On the estimation of ${\Cal Z}_2(s)$,
in ``Anal. Probab. Methods
Number Theory" (eds. A. Dubickas et al.), TEV, Vilnius, 2002,  83-98.

\item {[16]} A. Ivi\'c, On the integral of Hardy's function,  Arch. Mathematik
{\bf83}(2004), 41-47.

\item{[17]} A. Ivi\'c, The Mellin transform of the square of Riemann's
zeta-function, International J. of Number Theory {\bf1}(2005), 65-73.

\item {[18]} A. Ivi\'c, On some reasons for doubting the Riemann Hypothesis,
in P. Borwein, S. Choi, B. Rooney and A. Weirathmueller, ``The Riemann
Hypothesis'', CMS Books in Mathematics, Springer, 2008.

\item {[19]} A. Ivi\'c, On the Mellin transforms of powers of Hardy's function,
 Hardy-Ramanujan Journal {\bf33}(2010), 32-58.

\item {[20]} A. Ivi\'c, On some problems involving Hardy's function,
Central European J. Math. {\bf8(6)}(2010), 1029-1040.

\item {[21]} A. Ivi\'c, M. Jutila and Y. Motohashi, The Mellin
transform of powers of  the Riemann zeta-function,  Acta Arith.
{\bf95}(2000), 305-342.

\item{ [22]} A. Ivi\'c and Y. Motohashi,  A note on the mean value of
the zeta and L-functions VII, Proc. Japan Acad. Ser. A
{\bf 66}(1990), 150-152.

\item{ [23]} A. Ivi\'c and Y. Motohashi,  The mean square of the
error term for the fourth moment of the zeta-function, Proc. London Math.
Soc. (3){\bf 66}(1994), 309-329.

\item {[24]} A. Ivi\'c and Y. Motohashi,  The fourth moment of the
Riemann zeta-function,  J. Number Theory {\bf 51}(1995), 16-45.

\item{[25]}  M. Jutila, Mean values of Dirichlet series via Laplace
transforms,
in ``Analytic Number Theory" (ed. Y. Motohashi), London Math. Soc.
 LNS {\bf247},  Cambridge University Press, Cambridge, 1997, 169-207.

\item {[26]}  M. Jutila, The Mellin transform of the square of Riemann's
zeta-function,  Periodica Math. Hung. {\bf42}(2001), 179-190.

\item{[27]}  M. Jutila, The Mellin transform of the fourth power of the
Riemann zeta-function, in:
Adhikari, S.D. (ed.) et al., Number Theory. Proc.
Inter. Conf. on Analytic Number Theory with special emphasis
on $L$-functions, held at the Inst. Math. Sc., Chennai, India,
January 2002. Ramanujan Math. Soc. LNS {\bf1}2005), 15-29.

\item{[28]}  M. Jutila, Atkinson's formula for Hardy's function,
J. Number Theory {\bf129}(2009), 2853-2878.

\item{[29]} M. Jutila, An estimate for the Mellin transform of Hardy's function,
Hardy-Ramanujan J. {\bf33}(2010), 23-31.

\item{[30]} M. Jutila, The Mellin transform of Hardy's function is entire (in Russian),
Mat. Zametki {\bf88}(4)(2010), 635-639.

\item{[31]} M. Jutila, An asymptotic formula for the primitive of Hardy's
function,  in press in Arkiv Mat.  DOI:10.1007/s11512-010-0122-4

\item{[32]}  J.P. Keating  and  N.C.  Snaith, Random Matrix Theory
and $L$-functions at $s=1/2$,   Comm. Math. Phys. {\bf214}(2000), 57-89.

\item {[33]} H. Kober, Eine Mittelwertformel der Riemannschen Zetafunktion,
Compositio Math. {\bf3}(1936), 174-189.

\item{[34]}  M.A. Korolev, On the integral of Hardy's function $ Z(t)$,
 Izv. Math. {\bf72}, No. 3, 429-478 (2008); translation from  Izv.
Ross. Akad. Nauk, Ser. Mat. {\bf72}, No. 3, 19-68 (2008).

\item{ [35]} D.H. Lehmer,  On the roots of the Riemann zeta
function, Acta Math. {\bf 95}(1956), 291-298.

\item{ [36]} D.H. Lehmer, Extended computation of the Riemann
zeta-function, Mathematika {\bf 3}(1956), 102-108.

\item{[37]} M. Lukkarinen, The Mellin transform of the square of Riemann's
zeta-function and Atkinson's formula,  Doctoral Dissertation, Annales
Acad. Sci. Fennicae, No. {\bf140}, Helsinki, 2005, 74 pp.

\item{[38]} Y. Motohashi, Spectral  theory of the Riemann zeta-function,
Cambridge University Press, 1997.

\item{[39]} K.  Ramachandra, On the mean-value and omega-theorems
for the Riemann zeta-function, LN's {\bf85}, Tata Inst. of Fundamental Research
(distr. by Springer Verlag, Berlin etc.), Bombay, 1995.

\item{[40]} A. Selberg, Selected papers, Vol. 1,  Springer Verlag, Berlin, 1989.

\item{[41]}  E.C. Titchmarsh, Introduction to the theory of Fourier integrals
(2nd edition), Oxford University Press, Oxford, 1948.

\item{[42]}  E.C. Titchmarsh, The theory of the Riemann
zeta-function (2nd edition),  Oxford University Press, Oxford, 1986.

\item{[43]} N. Watt, A note on the mean square of $|\zt|$,
J. London Math. Soc. {\bf82(2)}(2010), 279-294.

\endRefs
\enddocument

\bye